\documentclass[preprint]{elsarticle}

\usepackage{amsmath}
\usepackage{graphics}

\begin{document}

\begin{frontmatter}

\title{%
The hypersecant Jacobian approximation for
quasi-Newton solves of sparse nonlinear systems}

\author{Johan Carlsson}
\ead{johan@txcorp.com}
\author{John R.~Cary}
\address{Tech-X Corporation, 5621 Arapahoe Avenue, Boulder, CO 80303, USA}

\begin{abstract}
A new Jacobian approximation is developed for use in
quasi-Newton methods for solving systems of nonlinear
equations. The new hypersecant Jacobian approximation is
intended for the special case where the evaluation of the
functions whose roots are sought dominates the computation time,
and additionally the Jacobian is sparse. One example of such
a case is the solution of the discretized transport equation
to calculate particle and energy fluxes in a fusion plasma.
The hypersecant approximation of the Jacobian is calculated
using function values from previous Newton iterations, similarly
to the Broyden approximation. Unlike Broyden, the hypersecant
Jacobian converges to the finite-difference approximation of
the Jacobian. The calculation of the hypersecant Jacobian
elements requires solving small, dense linear systems, where
the coefficient matrices can be ill-conditioned or even exactly
singular. Singular-value decomposition (SVD) is therefore used.
Convergence comparisons of the hypersecant method, the
standard Broyden method, and the colored finite differencing
of the PETSc SNES solver are presented.
\end{abstract}

\end{frontmatter}

\section{Introduction}

The Newton-Raphson method can provide solutions to systems
of nonlinear equations.  It works by finding a linear
approximation to the residual vector (the vector whose root
is desired), and then solving the resulting linear system to
obtain an approximate solution.  Iteration of this process,
if it converges, gives a numerical solution to the nonlinear
equations.  Thus, there are three parts to each iteration.
The first is to obtain the linear approximation, i.e., the
Jacobian, the matrix of the derivatives of the residual vector
components with respect to the parameters.  The second is
solving the resulting linear system.  The third is to correct
a change that makes the solution worse (e.g., by using a line
search along the direction of the change).  Here we concentrate
on the first two parts.

Generally, calculating the Jacobian is the most computationally
intensive part of the problem because the the Jacobian has $N^2$
terms for $N$ variables (and $N$ residuals).  In many cases of
interest it cannot be calculated directly, as the residuals are
not known as simple analytic functions.  Hence, the Jacobian is
found by the finite difference approximation.  However, when the
number of variables is large, this becomes prohibitive, hence
the development of quasi-Newton methods, such as the Broyden
method~\cite{broyden1965}. In the Broyden method, the Jacobian is
updated by adding a term to it so that it accurately predicts the
change in the residuals of the previous iteration.  However,
the Broyden method can be problematic; it is often found not
to converge, in large part due to this update not converging
to the true Jacobian.

Our motivating application is the time advance of systems
of nonlinear partial differential equations.  In such cases
there are a number of variables associated with each cell
upon discretization, and the residual is the change of those
variables in a given time step that gives one a stable advance.
Since the continuum equations are local, the residuals for a
cell depend on the values of the variables in that cell and
some nearby cells.  Hence, the Jacobian is sparse with $O(N)$
terms, and so it is reasonable to find the Jacobian through
finite differences. In particular, one can use \emph{coloring},
as is done in the PETSc SNES solver~\cite{petsc-efficient} to ensure
that the minimal number of evaluations are carried out to
obtain those finite differences.

We combine ideas from both these two approaches.
We recognize that the Jacobian is sparse, and so there are
only a small number of derivatives to be calculated for each
residual.  However, rather than explicitly computing all the
finite differences at each step, we use the results from the
last several solver iterations to solve for an approximation
to the Jacobian.  This may not be enough information for a
precise solution, as either there may be an insufficient number
of equations, or some of the equations may be degenerate (as
can happen when successive steps are in the same direction).
In the former case, which occurs for the first few iterations
of a quasi-Newton solve, we use Broyden until sufficient state
has been generated to calculate the hypersecant Jacobian.
In the latter case we use singular-value decomposition (SVD) to
determine the Jacobian.  Finally, because the resulting Jacobian
may in general itself be ill-conditioned, we solve for the linear
approximation for
the change in the parameters again through SVD.

We compare our method, the Broyden method, and the (colored)
PETSc SNES solver to three problems.  The first is a
three-variable linear system, which illustrates how the
Broyden Jacobian does not correctly converge.  In contrast,
our hypersecant method does provide the converged Jacobian.
For linear problems, once the Jacobian is correct,
the problem is solved in the next step.  Thus, for this case
there is no need to apply the PETSc SNES solver as the results
are known.  The second problem is a three-variable nonlinear
system.  There we find that the Broyden and hypersecant methods
are competitive, with the finite-difference approximation taking
more function evaluations.  This is perhaps not surprising, as
with only three variables, the Jacobian is nearly dense. Finally
we apply the methods to our motivating application, a nonlinear
transport equation of a kind that arises in certain problems
in plasma physics.  Here we find that the hypersecant method
takes roughly half the number of function evaluations to
converge compared with the finite difference approximation.
The Broyden convergence is much slower still.

The organization of this paper is as follows.  In the next
section we discuss the Broyden method for updating the Jacobian.
In the subsequent section we develop our hypersecant method.
Next we compare the methods in application to the simple
three-variable linear and nonlinear problems.  Subsequently
we compare them in application to our motivating problem,
a discretized transport equation.  Finally we summarize.

\section{Quasi-Newton methods and the Broyden approximation}

The goal of any Newton-like method is to solve a system of $N$
nonlinear equations
\begin{displaymath}
\mathbf{F}(\mathbf{x}) = 0
\end{displaymath}
of $N$ variables, $\mathbf{x}$. 
Given a current best guess $\mathbf{x}^k$
for the solution, one Taylor expands in its vicinity to obtain
the linear equation,
\begin{equation}
\mathbf{F}(\mathbf{x}^k + \delta\mathbf{x}^k) =
 \mathbf{F}(\mathbf{x}^k) + \mathbf{J}(\mathbf{x}^k) \cdot \delta\mathbf{x}^k +
  O(\delta\mathbf{x}^k \cdot \delta\mathbf{x}^k) = 0 \, ,
\label{eq:taylor}
\end{equation}
where $\mathbf{J}$ is the Jacobian
\begin{displaymath}
\mathbf{J}(\mathbf{x}) = \frac{\partial \mathbf{F}}{\partial \mathbf{x}} \, ,
\end{displaymath}
which implies that a better approximation $\mathbf{x}^{k+1}$ to the root
is found by solving for the linearized change $\delta\mathbf{x}^k$,
\begin{displaymath}
\mathbf{J}(\mathbf{x}^k) \cdot \delta\mathbf{x}^k = - \mathbf{F}(\mathbf{x}^k)
\, ,
\end{displaymath}
from which one obtains the improved guess for the root,
\begin{displaymath}
\mathbf{x}^{k+1} = \mathbf{x}^k + \delta\mathbf{x}^k
\, .
\end{displaymath}
Iterating this process then gives the root, if the process
converges.

In the standard Newton-Raphson method the true, analytic
Jacobian is used. If the true Jacobian is not known, a
finite-difference approximation can be used, making the scheme
a quasi-Newton method. However, if $\mathbf{F}(\mathbf{x})$
is numerically expensive to evaluate, the finite-difference
Jacobian can become inefficient due to the number of function
evaluations needed to compute it.

The Broyden method~\cite{broyden1965} is a common quasi-Newton
method where the Jacobian is approximated as:
\begin{equation}
\mathbf{J}^{k+1} = \mathbf{J}^{k} +
 \frac{\mathbf{F}^{k+1} - \mathbf{F}^{k} -
  \mathbf{J}^{k} \cdot
   \delta\mathbf{x}^k}{\delta\mathbf{x}^k \cdot \delta\mathbf{x}^k}
    \otimes \delta\mathbf{x}^k \, ,
\label{eq:broyden}
\end{equation}
where we have introduced the shorthand notation
$\mathbf{F}^k = \mathbf{F}(\mathbf{x}^k)$ and
$\mathbf{J}^k = \mathbf{J}(\mathbf{x}^k)$, and
where $\otimes$ is the dyadic (or outer) product.  The above
approximation comes from demanding that the new Jacobian
$\mathbf{J}^{k+1}$ predict the change in the residuals for
the most recent step,
\begin{equation}
\mathbf{F}^{k+1} - \mathbf{F}^k =
 \mathbf{J}^{k+1} \cdot (\mathbf{x}^{k+1} - \mathbf{x}^k)  \, .
\label{eq:broydreqsingle}
\end{equation}
The initial approximation $\mathbf{J}^0$ is typically set to
either the finite-difference Jacobian (good, but numerically
expensive) or the unit matrix (free).

The obvious advantage of the Broyden Jacobian is that it is
free in the sense that it requires no extra evaluations of
$\mathbf{F}$.  The disadvantage of the Broyden Jacobian is that
it is not guaranteed to converge toward the true Jacobian; in fact
error can accumulate until the increment $\delta\mathbf{x}$
found by inverting it is not even in the descent direction,
i.e. $\mathbf{x} + \delta\mathbf{x}$ is not closer to the root
than $\mathbf{x}$.  This occurs because Eq.~(\ref{eq:broydreqsingle}) is
underdetermined, so the update (\ref{eq:broyden}) is only one
of many possible.  In a solver using the Broyden approximation
one must intercept this occurrence and reinitialize the Jacobian.
Indeed, in our experience from using the Broyden method to solve
a discretized fusion transport equation~\cite{carlsson2002},
optimal performance, i.e. the minimal number of function
evaluations, required the Broyden Jacobian to be reinitialized
as often as every 5--10 Newton iterations.

\section{The hypersecant approximation of the Jacobian%
\label{sec:bhsalgo}}

Our hypersecant approximation is like the Broyden method in
that the Jacobian comes with no additional function evaluations,
but unlike Broyden, the approximate Jacobian will
approach the true Jacobian after sufficiently many iterations
as one approaches the solution, provided the directions of the
steps span the relevant space.  The primary difference is
that instead of requiring only \ref{eq:broydreqsingle}, we require
that this hold for several previous steps,
\begin{equation}
\mathbf{F}(\mathbf{x}^{k+1}) - \mathbf{F}(\mathbf{x}^{k-\ell}) =
 \mathbf{J}^{k+1} \cdot (\mathbf{x}^{k+1} - \mathbf{x}^{k-\ell})  \, , \
 \ell = 0, \ldots, L\!-\!1 \, .
\label{eq:broydreq}
\end{equation}

For the fully dense Jacobian, one obtains a full set of
equations after $N$ linear independent steps.  For large
systems, this is not practical in general.  However, it can be
very good when the Jacobian is sparse, such as occurs in the
implicit advance of discretized partial differential equations.
In particular, the above system of equations is determining
when the sparsity is such that any residual depends on only $L$
variables.  For example, in a finite difference discretization
of a one-dimensional, single-variable, second-order PDE,
each residual depends at most on three variables, and so after
three linearly independent steps one has a good approximation
to the Jacobian.

One issue is that the above approach cannot alone be used to
find a new Jacobian via Eq.~(\ref{eq:broydreq}) until one has
taken a number of steps equal to the maximum number of nonzero
elements in a row of the Jacobian. This is clearly the case for
the first few iterations. The Broyden update can then be used to
reduce the number of unknowns to make the system well determined,
see Section~\ref{sec:lineartest},until a sufficient number of steps
have been taken.
Even if a sufficient number of steps have been taken, the system
of equations~(\ref{eq:broydreq}) can be underdetermined, e.g.
because several successive steps are in nearly the same direction.  
This case is handled by using singular-value decomposition (SVD)
to solve for the hypersecant Jacobian elements, see
Section~\ref{sec:lineartest}.

\section{Comparison of hypersecant, Broyden and PETSc colored
FD Jacobian for few-variable systems%
\label{sec:comparison}}

In this section we compare the application of two approaches
(hypersecant and Broyden) to a linear three-variable system,
and then we compare the application of three approaches
(hypersecant, Broyden, and finite-difference Jacobian) to a
nonlinear three-variable system.  Our solvers did not attempt
to correct for an increment that takes the solution further
away from the true root.

\subsection{Linear system%
\label{sec:lineartest}}

The linear system is defined by the equations,
\begin{displaymath}
\left \{
\!\!
  \begin{array}{lcccccccccc}
    f_1(x_1, x_2 \quad \ \, )      & = &
            x_1     & + & \frac{x_2}{2} &   &               & - & \frac{3}{2} &
      = & 0 \\
    f_2(x_1, x_2, x_3) & = &
      \frac{x_1}{2} & + &       x_2     & + & \frac{x_3}{2} & - &          2  &
      = & 0 \\
    f_3(\quad \ \, x_2, x_3)   & = &
                    &   & \frac{x_2}{2} & + &       x_3     & - & \frac{3}{2} &
      = & 0 \\
  \end{array}
\right . \, .
\end{displaymath}
This system has the true solution $\mathbf{x} = (x_1, x_2, x_3) = (1, 1, 1)$,
and the Jacobian is
\begin{equation}
\mathbf{J} = \frac{d\mathbf{f}}{d\mathbf{x}} =
\left [
  \begin{array}{ccc}
    1   & 1/2 & 0   \\
    1/2 & 1   & 1/2 \\
    0   & 1/2 & 1   \\
  \end{array}
\right ] \, .
\label{eq:exactjac}
\end{equation}

A simple quasi-Newton solver was implemented using the
LAPACK~\cite{laug} implementation
ATLAS~\cite{WN147} version 3.8.3.
This simple linear system could of course have been solved in a more
direct manner, but the purpose of the test is primarily to compare
the Broyden and hypersecant Jacobian approximations.

We initialize the Jacobian as the unit matrix, $\mathbf{J}^0
= \mathbf{1}$ and use the initial guess $\mathbf{x} = (1/2,
1/2, 1/2)$. It takes 14 Newton iterations using the Broyden
approximate Jacobian to calculate the root to 12 significant
digits. The reason for the poor convergence is that the
Broyden Jacobian does not converge toward the true Jacobian.
At the final iteration, the Broyden Jacobian is
\begin{displaymath}
\mathbf{B}^{k=14} =
\left [
  \begin{array}{ccc}
    1.124366 & 0.599240 & 0        \\
    0.397825 & 0.823632 & 0.397825 \\
    0        & 0.599240 & 1.124366 \\
  \end{array}
\right ] \, ,
\end{displaymath}
which is far from the exact Jacobian~(\ref{eq:exactjac}).

For the hypersecant solve, we initialized the Jacobian
with the identity and used the Broyden update for the
first two iterations, but after the third iteration we used
Eqs.~(\ref{eq:broydreq}).  As an example, for the middle row, the
system of equations to solve is:
\begin{equation}
\left [
  \begin{array}{ccc}
    x_1^k\!-\!x_1^{k-3} \ &
    x_2^k\!-\!x_2^{k-3} \ &
    x_3^k\!-\!x_3^{k-3} \\
    x_1^k\!-\!x_1^{k-2} \ &
    x_2^k\!-\!x_2^{k-2} \ &
    x_3^k\!-\!x_3^{k-2} \\
    x_1^k\!-\!x_1^{k-1} \ &
    x_2^k\!-\!x_2^{k-1} \ &
    x_3^k\!-\!x_3^{k-1} \\
  \end{array}
\right ]
\left [
  \begin{array}{c}
    H_{2,1}^k \\
    H_{2,2}^k \\
    H_{2,3}^k \\
  \end{array}
\right ]
=
\left [
  \begin{array}{c}
    f_2^k\!-\!f_2^{k-3} \\
    f_2^k\!-\!f_2^{k-2} \\
    f_2^k\!-\!f_2^{k-1} \\
  \end{array}
\right ]
\, .
\end{equation}
Thus, at the third step we obtain the hypersecant Jacobian,
\begin{displaymath}
\mathbf{H}^{k=3} =
\left [
  \begin{array}{ccc}
    1.000000 & 0.500000 & 0        \\
    0.500000 & 1.000000 & 0.500000 \\
    0        & 0.500000 & 1.000000 \\
  \end{array}
\right ] \, ,
\end{displaymath}
It is seen to equal the finite-difference approximation of
the Jacobian (which in turn equals the true Jacobian for
the linear system used in this case). This simple example
thus clearly illustrates how the hypersecant Jacobian
avoids the major disadvantage of the Broyden approximation,
while maintaining the greatest advantage: the fact that it
is calculated without making any extra evaluations of the
function $\mathbf{f}(\mathbf{x})$.

This linear example also illustrates how singular-value
decomposition (SVD) is used to solve for the hypersecant
Jacobian elements. The variables $x_1$ and $x_3$ are
interchangeable. So if they are given the same initial values,
they will follow exactly the same trajectory. The linear system
to solve, Eqs.~(\ref{eq:broydreq}) for the hypersecant Jacobian
elements on row~1, therefore becomes singular:
{\scriptsize%
\begin{displaymath}
\begin{split}
&\left [
  \begin{array}{ccc}
    -4.143137616670 \times 10^{-2} &
    -1.429236794928 \times 10^{-1} &
    -4.143137616670 \times 10^{-2} \\
    -3.254780687737 \times 10^{-1} &
    -3.617952748235 \times 10^{-1} &
    -3.254780687737 \times 10^{-1} \\
    +4.245219312263 \times 10^{-1} &
    +6.382047251765 \times 10^{-1} &
    +4.245219312263 \times 10^{-1} \\
  \end{array}
\right ]
\times
\left [
  \begin{array}{c}
    H_{1,1}^3 \\
    H_{1,2}^3 \\
    H_{1,3}^3 \\
  \end{array}
\right ]
= \\
&\quad
\left [
  \begin{array}{c}
    -1.843550556595 \times 10^{-1} \\
    -6.872733435972 \times 10^{-1} \\
    +1.062726656403 \times 10^{+0} \\
  \end{array}
\right ]
\, .
\end{split}
\end{displaymath}%
}
Column~1 and 3 are identical and LAPACK \texttt{dgesvx()} fails
with error code \texttt{INFO = 4} signifying a coefficient
matrix singular to within machine precision. The inverse
condition number \texttt{RCOND = 7.426728596992e-17} confirms
this. The implementation of the hypersecant solver used here handles
situations like this by switching to SVD when \texttt{dgesvx()}
fails. LAPACK \texttt{dgesvd()} informs us that the singular
values of the coefficient matrix are
\begin{displaymath}
(w_1, w_2, w_3) =
 (1.060808064514, 9.516998001847 \times 10^{-2},
 4.162068222172 \times 10^{-17}) \, .
\end{displaymath}
Following the usual SVD prescription to solve linear systems,
our hypersecant solver implementation treats $w_3$ as
approximately zero and sets $1/w_3 = 0$ when it calculates
the hypersecant Jacobian elements
\begin{displaymath}
(H_{1,1}^3, H_{1,2}^3, H_{1,3}^3) =
 (0.500000000000, 1.000000000000, 0.5000000000000)
\end{displaymath}
which equals the true Jacobian to within machine
precision, exactly the answer we want.  The fact that
Eqs.~(\ref{eq:broydreq}) can be ill-conditioned or even exactly
singular therefore does not mean that their solutions, the
hypersecant Jacobian elements, are a poor approximation.
As shown here, by using SVD, hypersecant Jacobian elements
can be extracted, exact to within machine precision, even from
a very ill-conditioned system.

\subsection{Nonlinear system%
\label{seq:nonlineartest}}

The next test is for the system of nonlinear equations
\begin{equation}
\left \{
\!\!
  \begin{array}{lcccccccccc}
    f_1(x_1, x_2 \quad \ \, ) & = &
      \frac{x_1^2}{2} & + & \frac{x_2^2}{4} & & & - & \frac{3}{4} & = & 0 \\
    f_2(x_1, x_2, x_3) & = &
      \frac{x_1^2}{4} & + & \frac{x_2^2}{2} & + & \frac{x_3^2}{4} & - & 1  &
      = & 0 \\
    f_3(\quad \ \, x_2, x_3) & = &
      & & \frac{x_2^2}{4} & + & \frac{x_3^2}{2} & - & \frac{3}{4} & = & 0 \\
  \end{array}
\right . \, .
\label{eq:snesbhs2}
\end{equation}
This system also has the root
$\mathbf{x} = (x_1, x_2, x_3) = (1, 1, 1)$, and the Jacobian at the root is
\begin{displaymath}
\mathbf{J}(1, 1, 1) =
\left [
  \begin{array}{ccc}
    1   & 1/2 & 0   \\
    1/2 & 1   & 1/2 \\
    0   & 1/2 & 1   \\
  \end{array}
\right ] \, .
\end{displaymath}
To be able to compare not only the hypersecant and Broyden
Jacobians, but also a finite-difference (FD) Jacobian,
a quasi-Newton solver was implemented using the PETSc
SNES nonlinear solver toolkit. PETSc can do colored finite
differencing, i.e. take full advantage of the sparsity of the
Jacobian to minimize the number of function evaluations used
to calculate the FD Jacobian. We only specify the sparsity
pattern of the Jacobian and let PETSc determine how to do the
colored finite differencing.

For the hypersecant Jacobian, we take advantage of the fact
that the first and last rows of the Jacobian have only two
non-zero elements. For the first row, the hypersecant Jacobian
elements are given as the solution to the linear system
\begin{equation}
\left [
  \begin{array}{cc}
    x_1^k\!-\!x_1^{k-2} \ &
    x_2^k\!-\!x_2^{k-2} \ \\
    x_1^k\!-\!x_1^{k-1} \ &
    x_2^k\!-\!x_2^{k-1} \ \\
  \end{array}
\right ]
\left [
  \begin{array}{c}
    H_{1,1}^k \\
    H_{1,2}^k \\
  \end{array}
\right ]
=
\left [
  \begin{array}{c}
    f_1^k\!-\!f_1^{k-2} \\
    f_1^k\!-\!f_1^{k-1} \\
  \end{array}
\right ]
\, .
\end{equation}
For the first and the last rows of the Jacobian, the
hypersecant approximation can thus be used already after two
iterations (the hypersecant Jacobian elements can be calculated
from the available state at the end of the second iteration,
and is thus first available for use during the third iteration).

When the iteration number is too small to calculate the
hypersecant Jacobian elements ($k=1$ for row 1 and 3, and
$k=1,2$ for row 2), we first do a Broyden update of the old
Jacobian. Elements of this updated Jacobian are then used to
replace sufficiently many unknowns in Eq.~(\ref{eq:broydreq}) to
reduce the number of unknowns to equal the number of iterations,
thereby making the system well determined. For example, for
the first iteration we get for the first row of the Jacobian
\begin{displaymath}
\tilde{H}_{1,1}^1 (x_1^1 - x_1^0) + B_{1,2}^1 (x_2^1 - x_2^0) =
 f_1^1 - f_1^0 \, ,
\end{displaymath}
which only has the single unknown $\tilde{H}_{1,1}^1$, where
the tilde is used to indicate that this is only an approximate
hypersecant Jacobian element. This equation can be trivially
solved,
\begin{displaymath}
\tilde{H}_{1,1}^1 =
\frac{f_1^1 - f_1^0 - B_{1,2}^1 (x_2^1 - x_2^0)}%
{(x_1^1 - x_1^0)} \, .
\end{displaymath}
Because the diagonal elements have larger values, we
somewhat arbitrarily choose to use $B_{1,2}^1$ to calculate
$\tilde{H}_{1,1}^1$, instead of $B_{1,1}^1$ to calculate
$\tilde{H}_{1,2}^1$.  Similarly, for the middle row we get
\begin{displaymath}
\tilde{H}_{2,2}^1 =
\frac{f_2^1 - f_2^0 - B_{2,1}^1 (x_1^1 - x_1^0) - B_{2,3}^1 (x_3^1 - x_3^0)}%
{(x_2^1 - x_2^0)} \, ,
\end{displaymath}
etc. For $k \geq 3$, the full hypersecant Jacobian is
calculated (and subsequently used for $k \geq 4$).

\begin{figure}
\begin{center}
\includegraphics[width=.9\columnwidth]{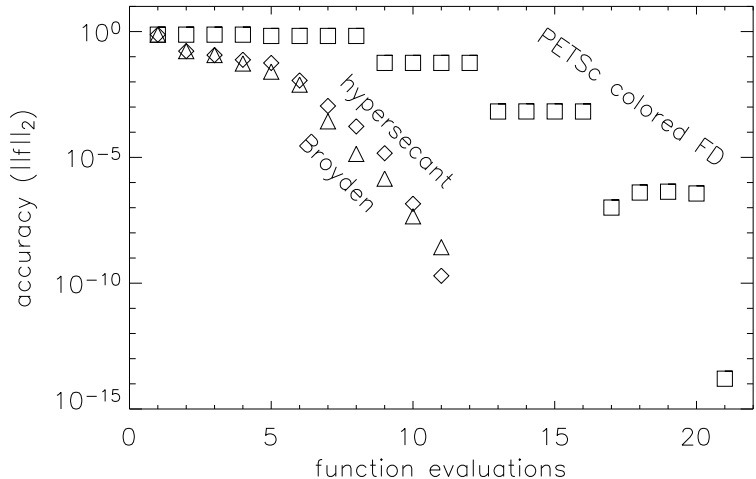}
\end{center}
\caption{Accuracy as a function of the number of
function evaluations when doing quasi-Newton solve of
Eq.~(\ref{eq:snesbhs2}).  For the triangles, the Broyden
approximation of the Jacobian was used.  For the diamonds,
the hypersecant Jacobian. For the squares, the PETSc colored
finite-difference Jacobian.}
\label{fig:simple}
\end{figure}

As usual, we initialize the Jacobian as the unit matrix and
use the initial guess $\mathbf{x} = (1/2, 1/2, 3/2)$. As
can be seen in Fig.~\ref{fig:simple}, the quasi-Newton solve
converges to the default SNES accuracy in eleven iterations when
either the hypersecant or Broyden Jacobians are used. Using
the SNES colored FD approximation of the Jacobian, the solve
converges in just five iterations. However, because the colored
FD approximation requires three extra function evaluations
(for a Jacobian with this particular sparsity pattern), the
total number of function evaluations is 21 vs. the 11 needed
by hypersecant/Broyden.

Even if the final Broyden Jacobian is a much poorer
approximation of the true Jacobian than the final hypersecant
Jacobian is, the performance of the Broyden and hypersecant
methods in this case is very similar. However, the hypersecant
Jacobian gives more rapid convergence for the last
few iterations, as it becomes a more accurate approximation
to the true Jacobian.

%
%

\section{Application to a nonlinear transport equation%
\label{seq:transporttest}}

As noted in the introduction, our motivation for this research is
to develop methods for solving nonlinear transport equations.
Our particular example comes from cross-flux-surface transport in
a fusion confinement device, where the transport equation is a
radial continuity equation~\cite{hirshman1979}:
\begin{equation}
\frac{\partial \mathbf{u}}{\partial t} +
\frac{1}{V'(r)} \frac{\partial}{\partial r} V'(r) \, \pmb{\Gamma} =
 \mathbf{S} \, ,
\label{eq:xport}
\end{equation}
where the highly nonlinear flux is given by $\pmb{\Gamma} =
\pmb{\Gamma}(r, \mathbf{u}, \mathbf{u}')$ and $\mathbf{u}' =
\partial \mathbf{u} / \partial r$. The volume element, the
differential volume in the radial variable $r$, is denoted by
$V'(r)$.  The components of $\mathbf{u}(r)$ are typically the
temperatures and densities of the plasma particle species, but
they could also include momenta. In the cases where it is
difficult to advance this partial differential equation, the
dominant contribution to the flux comes from turbulent transport,
which is triggered when a critical gradient is exceeded. This
creates a sensitive dependence on $\mathbf{u}'$. The stiffness of
the transport equation makes it necessary to use implicit time
discretization.


Time-implicit, spatially second-order finite-difference
discretization on an equidistant mesh $r_j = j\Delta r, j = 0,
\ldots, N$ gives the discretized system of nonlinear equations:
\begin{equation}
\begin{split}
&\mathbf{F}_j(\Delta\mathbf{u}_{j-1}, \Delta\mathbf{u}_{j},
 \Delta\mathbf{u}_{j+1}) =
 \Delta\mathbf{u}_{j} + \\
&\theta \left ( \frac{1}{W'(r_j)}
 \frac{W'(r_{j+1/2}) \pmb{\Gamma}^{n+1}_{j+1/2} -
       w'(r_{j-1/2}) \pmb{\Gamma}^{n+1}_{j-1/2}}{\Delta r} -
  \mathbf{S}^{n+1}_j
 \right ) \Delta t + \\
&(1 - \theta) \left ( \frac{1}{W'(r_j)}
 \frac{W'(r_{j+1/2}) \pmb{\Gamma}^{n}_{j+1/2} -
       W'(r_{j-1/2}) \pmb{\Gamma}^{n}_{j-1/2}}{\Delta r} -
  \mathbf{S}^{n}_j
 \right ) \Delta t = 0 \, , \\
\end{split}
\label{eq:discretized}
\end{equation}
where subscript denotes spatial location and superscript time
level, e.g. $\mathbf{S}^n_j = \mathbf{S}(r_j, n\Delta_t)$ is the
source at radius $r = r_j$ and time $t = n\Delta t$. The
independent variables $\Delta\mathbf{u}_j = \mathbf{u}^{n+1}_j -
\mathbf{u}^n_j$ and $\pmb{\Gamma}_{j+1/2} =
\pmb{\Gamma}(r_{j+1/2}, \mathbf{u}_{j+1/2},
\mathbf{u}'_{j+1/2})$, with $\mathbf{u}_{j+1/2} = (\mathbf{u}_{j}
+ \mathbf{u}_{j+1})/2$ and $\mathbf{u}'_{j+1/2} =
(\mathbf{u}_{j+1} - \mathbf{u}_{j})/\Delta r$.  The
implicitness parameter $\theta = 0$ gives fully explicit (FE)
equations, $\theta = 1$ fully implicit (FI) and $\theta = 1/2$
time-centered implicit (CI), with second-order temporal accuracy.

We next compare the three Jacobian approximations for this
transport equation~(\ref{eq:discretized}) with a single profile
evolved [e.g. the ion temperature $T_i(r)$].  For the nonlinear
flux we use the linear critical gradient model, which assumes the
flux $\Gamma = - \chi u'$ and diffusivity
\begin{displaymath}
\chi = \max\{(1/L - 1/L_c)/L, \chi_{min}\} \, ,
\end{displaymath}
where $L = u / |u'|$ and $L_c$ and $\chi_{min}$ are parameters.
We also use a source term $S(r) = 1 - r^\alpha$. We somewhat arbitrarily
choose the parameters $L = 1/2$, $\chi_{min} = 1/10$ and $\alpha = 2$.
With fully implicit time discretization we get the system of nonlinear
equations
\begin{equation}
\begin{split}
F_j(&\Delta u_{j-1}, \Delta u_{j}, \Delta u_{j+1}) = \Delta u_{j} + \\
&\left ( \frac{1}{W'(r_j)}
\frac{W'(r_{j+1/2}) \Gamma^{n+1}_{j+1/2} - W'(r_{j-1/2}) \Gamma^{n+1}_{j-1/2}}%
{\Delta r} - S^{n+1}_j \right ) \Delta_t = 0 \, , \\
\end{split}
\label{eq:snesbhslcgm}
\end{equation}
$j = 1, \ldots, N-1 \, .$
To determine the appropriate on-axis ($r_0 = 0$) boundary condition,
we take one step back and write the flux term as an undiscretized
divergence term:
\begin{displaymath}
W'(r_0) F_0 =  W'(r_0) \Delta u_0 +
\left ( \frac{\partial}{\partial r} W'(r) \Gamma \right )_0 \Delta t -
 W'(r_0) S^{n+1}_0 \Delta t = 0 \, ,
\end{displaymath}
where we have also multiplied the equation with $W'(r_0)$. In the limit
$r_0 \rightarrow 0$, in which $W'(r_0) \rightarrow 0$, we get
\begin{displaymath}
\left ( \frac{\partial}{\partial r} W'(r) \Gamma \right )_0 \Delta t = 0 \, .
\end{displaymath}
Discretizing this equation with second-order error gives
\begin{equation}
3 \Gamma^{n+1}_{1/2} - \Gamma^{n+1}_{3/2} = 0 \, ,
\end{equation}
which is a nonlinear on-axis boundary condition. The first term
depends on $\Delta u_0$ and $\Delta u_1$, and the second term on
$\Delta u_1$ and $\Delta u_2$. The first row of the Jacobian will
thus have non-zero elements in the first three columns. The
Jacobian is then tri-diagonal, but with $J_{0,2} \neq 0$, due to
the nonlinear BC.

The hypersecant Jacobian is calculated in the way described above in
Section~\ref{seq:nonlineartest}, but both hypersecant and Broyden
Jacobians are initialized as the unit matrix with the first row replaced
by $J_{0,0} = 3/10 \times N$, $J_{0,1} = -4/10 \times N$ and
$J_{0,2} = 1/10 \times N$ (the values given by the nonlinear BC
at steady state).

\begin{figure}
\begin{center}
\includegraphics[width=.9\columnwidth]{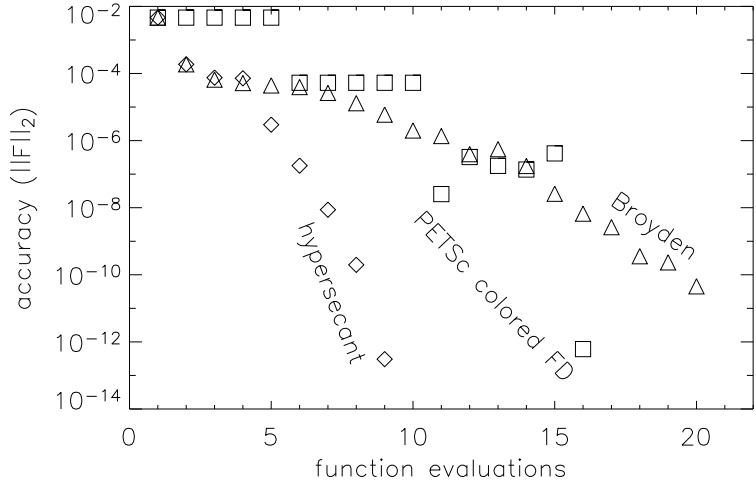}
\end{center}
\label{fig:lcgm}
\caption{Accuracy as a function of the number of function evaluations
when doing quasi-Newton solve of Eq.~(\ref{eq:snesbhslcgm}).
For the triangles, the Broyden approximation of the Jacobian was used.
For the diamonds, the hypersecant Jacobian. For the squares, the PETSc
colored finite-difference Jacobian.}
\end{figure}

Taking a single time step of 0.1~ms, the PETSc SNES quasi-Newton solve
converges as shown in Fig.~\ref{fig:lcgm} for the three different
Jacobian approximations.
%
%
With the hypersecant Jacobian, the solve converges using 9
function evaluations.  With PETSc colored FD 16, and with Broyden
20 function evaluations.  With PETSc colored FD, only three
quasi-Newton iterations are needed, but during each iteration,
four extra functions calls are used to calculate the colored-FD
Jacobian. With hypersecant Jacobian, the convergence is poor
until the fourth iteration, the first iteration where the full
hypersecant Jacobian is used. After that it takes five more
iterations until convergence (for a total of nine), but each
iteration only requires a single function evaluation. With the
Broyden Jacobian the convergence is much poorer.


For larger time steps, a problem with using the hypersecant
Jacobian approximation for PETSc SNES quasi-Newton solves becomes
evident.  The hypersecant Jacobian is a poor approximation of
the true Jacobian for the first few iterations. To use
hypersecant quasi-Newton solves for the transport equation with
large time steps, one probably even needs to accept that the
accuracy gets worse for the first few iteration.  The purpose of
the first few iterations is not primarily to get closer to the
root, but to sample enough function values to make
Eq.~(\ref{eq:broydreq}) well determined to allow the calculation
of the full hypersecant Jacobian. PETSc SNES tries to improve
the accuracy for every iteration by doing a line search for the
minimum accuracy along the direction of the increment. This
typically fails for a poor approximation of the Jacobian. We are
currently investigating some promising ways to solve these
problems and plan to report on our findings in a future
publication.

\section{Conclusions}

We have introduced the hypersecant Jacobian approximation, which
has improved characteristics for solving systems of nonlinear
equations with sparse Jacobians, such as occurs in the implicit
advance of nonlinear PDEs.  The basic idea behind this
approximation is to keep several historical values -- sufficient
to be able to solve uniquely for an approximate Jacobian
generically, but with use of Singular Value Decomposition to find
a best solution when the equations are degenerate.  Like the Broyden
Jacobian approximation, the hypersecant Jacobian is calculated without
making any additional function evaluations. Unlike the Broyden Jacobian
it converges toward the finite-difference approximation of the Jacobian.

Multiple examples (linear and nonlinear systems of three
independent variables and the 1D transport equation with
linear-critical-gradient diffusivity) were solved using these
different methods.  The new hypersecant method was shown to lead to
a convergent Jacobian.  Improved overall convergence with
respect to number of function evaluations was observed.

\vspace{4.0ex}
\noindent%
\textbf{Acknowledgement:}
Work supported by the United States Department of Energy
Office of Fusion Energy Sciences Grant \#DE-FG02-05ER84383.


\end{document}